\documentclass{amsart}

\usepackage{geometry} 
\usepackage{amssymb}
\usepackage{relsize}
\usepackage[all]{xy}
\usepackage{xcolor}
\usepackage{dsfont}
\usepackage[hidelinks]{hyperref}
\usepackage{cleveref}

\numberwithin{equation}{section}

\newtheorem{theorem}[subsubsection]{Theorem}
\newtheorem{proposition}[subsubsection]{Proposition}

\newtheorem{lemma}[subsubsection]{Lemma}

\newtheorem*{theorem*}{Theorem}

\theoremstyle{definition}
\newtheorem{definition}[subsubsection]{Definition}

\theoremstyle{remark}

\newtheorem{remark}[subsubsection]{Remark}

\title[Naturality of cohomological support]{On naturality of cohomological support for small quantum groups}
\author{Kevin Coulembier}
\address{School of Mathematics and Statistics, University of Sydney, Australia}
\email{kevin.coulembier@sydney.edu.au}

\newcommand{\Sym}{\operatorname{Sym}}

\newcommand{\bk}{\Bbbk}
\newcommand{\cN}{\mathcal{N}}

\newcommand{\Tilt}{\mathtt{Tilt}}

\newcommand{\Res}{\operatorname{Res}}

\newcommand{\mZ}{\mathbb{Z}}
\newcommand{\mN}{\mathbb{N}}

\newcommand{\mG}{\mathbb{G}}

\newcommand{\hyp}{\mathrm{hyp}}

\newcommand{\pd}{\mathrm{pd}}

\newcommand{\Hom}{\operatorname{Hom}}

\newcommand{\tto}{\twoheadrightarrow}
\newcommand{\Ext}{\operatorname{Ext}}

\newcommand{\Supp}{\operatorname{Supp}}
\newcommand{\mProj}{\operatorname{mProj}}
\newcommand{\Proj}{\operatorname{Proj}}

\newcommand{\fm}{\mathbf{m}}
\newcommand{\fp}{\mathbf{p}}
\newcommand{\fq}{\mathbf{q}}

\newcommand{\ffp}{\mathfrak{p}}

\newcommand{\Spec}{\operatorname{Spec}}

\newcommand{\cV}{\mathcal{V}}
\newcommand{\mD}{\mathbb{D}}
\newcommand{\mP}{\mathbb{P}}

\newcommand{\fg}{\mathfrak{g}}
\newcommand{\fb}{\mathfrak{b}}
\newcommand{\fl}{\mathfrak{l}}
\newcommand{\fn}{\mathfrak{n}}
\newcommand{\mC}{\mathbb{C}}

\newcommand{\supp}{\operatorname{supp}}

\newcommand{\DK}{\mathrm{DK}}

\renewcommand{\L}{\mathrm{L}}

\newcommand{\Mat}{\operatorname{Mat}}
\begin{document}

\begin{abstract}
We apply the theory of hypersurface support by Negron and Pevtsova to prove the naturality of cohomological support of the small quantum Borel groups, with respect to restriction to Levi subalgebras. This in turn implies naturality of cohomological support of the small quantum groups for modules that lift to Lusztig's big quantum group, for instance tilting modules. This allows us to make the recent conditional classification of tensor ideals of abelian type in the category of such tilting modules unconditional.
\end{abstract}


\subjclass[2020]{}

\maketitle

\section*{Introduction}

A powerful result in the theory of cohomological support varieties for restricted Lie algebras (in positive characteristic) connects it to a rank variety approach on the restricted nullcone, see \cite[\S 1]{FP}. Despite the profound similarities with cohomological support for small quantum groups, see \cite{GK, OsSupp, BNPP, Dr1}, no obvious such connection is available in the latter situation. 

As a consequence, the clean connection between cohomological support of a module and that of its restriction to an appropriate subalgebra does not follow immediately for small quantum groups the way it does for restricted Lie algebras. Recently this presented a direct obstacle in a classification effort concerning tilting modules for Lusztig's form of the big quantum group in \cite{CEO}.

More concretely, given a semisimple complex Lie algebra $\fg$ and a complex root $\zeta$ of $1$, satisfying standard compatibility conditions as made precise in the main text, Lusztig's small quantum group $u_\zeta(\fg)$ has cohomology ring $H^\bullet(u_\zeta(\fg),\mC)\cong\mC[\cN]$, with $\cN\subset\fg$ the nilpotent cone. Any finite dimensional $u_\zeta(\fg)$-module has cohomological support variety $\supp M$ in $\cN$. For a (standard) Levi subalgebra $\fl\subset\fg$, the support of the restriction to $\fl$ is always contained in $\supp M\cap \fl$, and one actually expects equality. Such equality is referred to as `naturality' of the support theory.

Rather than relying on a rank variety approach\footnote{An oversimplification of the connection between the ideas of the rank varieties and hypersurface support for a Hopf algebra $H$ is as follows. For the former one tries to understand whether a point is in the support by considering whether the module is projective (= has finite projective dimension) over a subalgebra of $H$ generated by a corresponding nilpotent element. For the latter, one tries to `resolve all other nilpotents' in $H$ and one considers whether the module has finite projective dimension over the resulting bigger algebra that surjects onto $H$.}, we apply ideas of `hypersurface support' due to Negron and Pevtsova \cite{NP1, NP2}. We show that this theory lends itself to naturality in Proposition~\ref{prop:iteratedOre2}. By results from \cite{NP2}, cohomological support for small quantum groups of Borel subalgebras $\fb\subset\fg$ of semisimple Lie algebras can be expressed in terms of hypersurface support. Combining this idea with our Proposition~\ref{prop:iteratedOre2} then shows that the cohomological support of $u_\zeta(\fb)$-modules restricted to $u_\zeta(\fb\cap\fl)$ exhibits naturality.

If we are only interested in the cohomological support, with respect to $u_\zeta(\fg)$, of modules over Lusztig's big quantum group $U_\zeta(\fg)$, such as tilting modules, then using ideas from \cite{BKN} we can prove the desired naturality by a reduction to the Borel subalgebra. The main observation is that the support in $\cN$ is closed under action of the adjoint group, so that it is expressible in terms of nilpotent orbits, which always intersect the Borel subalgebra non-trivially.

The paper is organised as follows. In Section~\ref{sec:ore} we discuss Ore extensions and present some results that will be used later when we view the De Concini-Kac forms of the quantum Borel groups as iterated Ore extensions. In Section~\ref{sec:NP} we present the main ideas of hypersurface support from \cite{NP1, NP2} that we will rely on, and connect them to Ore extensions. In Section~\ref{sec:Borel}, we derive naturality for the small quantum Borel groups. In Section~\ref{sec:Full} we derive naturality for the small quantum groups for modules that lift to Lusztig's big quantum group. Finally, in Appendix~\ref{sec:app} we justify our rephrasing of the results from \cite{NP1, NP2} by establishing equivalence between our polynomial approach and the formal approach from \cite{NP1, NP2}.

\subsection*{Use of AI} The fact that the theory from \cite{NP1, NP2} could be useful for proving naturality of cohomological support for small quantum Borel groups (Theorem~\ref{thm:B}) was discovered during experimentation in small rank cases with GPT-6 Astra. The final write-up is entirely due to the author who takes full responsibility for the text and validity of the results. 


\section{Ore extensions}\label{sec:ore}
Let $K$ be a commutative ring.

\subsection{Definitions}
We recall the notion of an Ore extension, see for instance \cite{GW}.

\begin{definition}
A morphism $R\to S$ of $K$-algebras is an {\bf Ore extension} if there exists $x\in S$ such that $S$ is free as a left and as a right $R$-module with basis $\{x^j\mid j\in\mN\}$, and such that
\[xR\subset Rx\oplus R\quad\mbox{and}\quad Rx\subset xR\oplus R.\]
\end{definition}

Using the left-hand side in the display above, we find $K$-linear maps $\sigma:R\to R$ and $\delta:R\to R$ defined via
\[xr=\sigma(r)x+\delta(r),\]
for all $r\in R$, where $\sigma$ is a $K$-algebra automorphism and $\delta$ is a $\sigma$-derivation:
\[\delta(r_1r_2)=\sigma(r_1)\delta(r_2)+\delta(r_1)r_2,\qquad\mbox{for all}\; r_1,r_2\in R.\]

Starting from a $K$-algebra automorphism $\sigma:R\to R$ and a $\sigma$-derivation $\delta:R\to R$ for a $K$-algebra $R$, there always exists a corresponding Ore extension $S=R[x;\sigma,\delta]$.

\begin{definition}
A $K$-algebra morphism $\phi:R\to S$ is an {\bf iterated Ore extension} if there exists $n\in\mN$ such that $\phi$ is a composite morphism
\[R=R_0\to R_1\to R_2\to\cdots\to R_n=S,\]
where each $R_j\to R_{j+1}$ is an Ore extension.
\end{definition}
If we want to be more precise, we can say that the above is an $n$-iterated Ore extension, so that $1$-iterated Ore extensions are precisely Ore extensions.
The following lemma is immediate, but it will be useful to have it written out.
\begin{lemma}\label{lem:Orequo}
Let $R\to S$ be an iterated Ore extension and consider $z\in R$ which is central in $S$ (and so also in $R$). Then $R/(z)\to S/(z)$ is an iterated Ore extension.
\end{lemma}

\subsection{Projective dimension}

\begin{proposition}\label{prop:iteratedOre}
Let $R\to S$ be an iterated Ore extension. For any left $S$-module $M$, we have
$\pd_RM<\infty$ if and only if $\pd_SM<\infty$. Moreover, if $R\to S$ is an $n$-iterated Ore extension, then the global dimension of $S$ is bounded by the global dimension of $R$ plus $n$.
\end{proposition}

By iteration, it suffices to deal with the case of an Ore extension. The proposition thus follows from the following lemma.
\begin{lemma}\label{lem:Ore}
Let $R\to S$ be an Ore extension. For any left $S$-module $M$, we have
\[\pd_RM\le\pd_SM\le \pd_R M+1.\]
\end{lemma}
\begin{proof}
By assumption, $S$ is free as a left $R$-module, so that any projective $S$-module remains projective as an $R$-module, proving the left inequality. To prove the right inequality, we choose a presentation $S=R[x;\sigma,\delta]$.

Denote by $\widetilde{M}$ the $R$ module obtained from $M$ that has the same underlying $K$-module but for which the new action of $r\in R$ is given by the old action of $\sigma^{-1}(r)$. Then we have a short exact sequence of $S$-modules
\[0\to S\otimes_R\widetilde{M}\xrightarrow{d} S\otimes_RM\xrightarrow{p} M\to 0,\]
where $d(s\otimes v)=sx\otimes v-s\otimes xv$ and $p(s\otimes v)=sv$. For an arbitrary $S$-module $N$ and $r\in\mN$ we thus find an exact sequence
\[\Ext^r_{R}(\widetilde{M},N)\to \Ext^{r+1}_S(M,N)\to \Ext^{r+1}_R(M,N).\]
The second inequality now follows from $\pd_R\widetilde{M}=\pd_RM$.
\end{proof}


\section{Our support group}\label{sec:NP}

All notions in this section can be extended to a scheme-theoretic set-up. Since that will not be required for our applications, we only report on simplified versions of the work in \cite{NP1,NP2}. We work over an algebraically closed field $\bk$.

\subsection{Cohomological support}

\subsubsection{Ring-theoretic support}
Let $S$ be a commutative non-negatively graded $\bk$-algebra  with $S_0=\bk$. For a graded $S$-module $M$, we set
\[\Supp_S(M)\;=\;\{\fm\in\mProj S\mid M_{\fm}\not=0\},\]
where $\mProj$ refers to the set of closed points in $\Proj S$. For a graded $S$-algebra $\phi:S\to R$ viewed as an $S$-module, we have
\begin{equation}\label{eq:suppphi}\Supp_S(R)\;=\; V(\ker \phi),\end{equation}
where $V(I)\subset\mProj S$ stands for the closed subset defined by a homogeneous ideal $I$.

\subsubsection{Cohomological support}
For a Hopf algebra $A$ over $\bk$, consider the (commutative) ext-algebra
\[E_A:=\Ext^{2\bullet}_A(\bk,\bk)=H^{2\bullet}(A,\bk).\]
Here $\bk$ will always stand for the trivial representation of a Hopf algebra, given by the co-unit $\epsilon:A\to\bk$.
We always view $E_A$ as a graded subalgebra of the full extension algebra, so by construction it is contained in even degrees.
For any $A$-module $M$, its cohomological support is given by
\[\cV_A(M)\;:=\;\Supp_{E_A}(\Ext^\bullet_A(M,M))\;\subset\; \mProj E_A,\]
where the module structure comes from the algebra morphism $E_A\to\Ext^\bullet_A(M,M)$, $\eta\mapsto  \eta\otimes M$, which is compatible with for instance \cite{Dr1}, but differs from the convention $M\otimes \eta$ in \cite{NP2}.

\begin{lemma}\label{lem:Hopfnormal}
Let $H$ be a finite dimensional Hopf algebra, with a normal Hopf subalgebra $A$, such that the quotient Hopf algebra $B:=H/\hspace{-1mm}/A= H/(A_+)$ is semisimple and acts trivially on $E_{A}$, then $E_A=E_H$. Moreover, as subsets of $\mProj E_A=\mProj E_H$, for an $H$-module $M$, we have
\[\cV_A(\Res^H_AM)\;=\; \cV_H(M).\]
\end{lemma}
\begin{proof}
By the collapsing Lyndon–Hochschild–Serre spectral sequence, for any $H$-module $N$ we have the restriction isomorphism
\[\Ext^\bullet_H(N,N)=\Hom_{B}(\bk,\Ext^\bullet_A(N,N)).\]
Applying this to $N=\bk$ shows $E_A=E_H$. Furthermore, the ring homomorphism $E_A\to\Ext^\bullet_A(M,M)$ factors via the inclusion $\Ext^\bullet_H(M,M)\hookrightarrow \Ext^\bullet_A(M,M)$, so that the conclusion follows from \eqref{eq:suppphi}.
\end{proof}

\subsection{Hypersurface support}\label{sec:hyper}

\subsubsection{} For a positive integer $n$, we set $Z_n=\bk[y_1,y_2,\ldots, y_n]$, viewed as an augmented algebra via $Z\to\bk, y_i\mapsto 0$ with kernel $\fm_n=(y_1,\ldots, y_n)$. We consider the $n$-dimensional space $W_n:=\fm_n/\fm_n^2$.

 When clear from context, we will leave out the subscripts $n$.

\subsubsection{}\label{def:datum} For a positive integer $n$, we define an \emph{$n$-deformation datum} $\mD$ to be a central $\bk$-algebra morphism $Z=Z_n=\bk[y_1,\ldots, y_n]\to Q$, where $Q$ is finite free over $Z$ (consequently, $Q$ is noetherian) and $Q$ has finite global dimension. This datum in particular defines a finite dimensional $\bk$-algebra $R:=\bk\otimes_{Z}Q\cong Q/\fm Q$.

\subsubsection{}\label{def:hyper} Consider an $n$-deformation datum $\mD$ with $Q\tto R$ as above. For a finite dimensional $R$-module $M$, following \cite{NP1, NP2} we define
\[\Supp^{\hyp}_{\mD}(M)\;=\;\{x\in\mP W\mid \pd_{Q/(f_x)}M=\infty\}\subset\mP W,\]
where $f_x\in\fm\subset Q$ stands for an arbitrary element that maps to $x$ under $\fm\tto\fm/\fm^2=W\dashrightarrow\mP W$. The fact that the definition does not depend on the choice of $f_x$ follows from  \cite[Corollary~5.5]{NP1} or \cite[Corollary~5.3]{NP2}.

\begin{remark}
Strictly speaking, \cite{NP1,NP2} do not apply to our setting, but instead to the $\fm$-adic completions $\widehat{Z}$ and $\widehat{Q}$. We explain in Appendix~\ref{sec:app} that this is not a genuine issue for the results of \cite{NP1, NP2} that we state in this section.
\end{remark}

\subsubsection{}\label{sec:DD}
We consider the following set-up, for positive integers $n_1\le n_2$. Consider a commutative diagram of algebra morphisms
\begin{equation}\label{eq:DD}\xymatrix{
Z_{n_1}\ar[rr]\ar[d]&& Q_1\ar[d]\\
Z_{n_2}\ar[rr]&& Q_2,
}\end{equation}
where the horizontal arrows are deformation data $\mD_1,\mD_2$ and the left vertical map is the canonical inclusion $y_i\mapsto y_i$. The latter induces a $\bk$-linear inclusion $W_{n_1}\to W_{n_2}$. We also obtain an induced algebra morphism $R_1\to R_2$.

The following proposition shows why hypersurface support is a useful tool to establish naturality.
\begin{proposition}\label{prop:iteratedOre2} Assume that $Q_1\to Q_2$ in \eqref{eq:DD} is an iterated Ore extension.
For a finite dimensional $R_2$-module $M$, we have
\[\Supp^{\hyp}_{\mD_1}(\Res^{R_2}_{R_1} M)\;=\;\Supp^{\hyp}_{\mD_2}(M)\cap \mP W_{n_1}.\]
\end{proposition}
\begin{proof}
Take $x\in \mP W_{n_1}\subset\mP W_{n_2}$ and consider a corresponding lift $f_x\in \fm_{n_1}\subset\fm_{n_2}$. By Lemma~\ref{lem:Orequo}, $Q_1/(f_x)\to Q_2/(f_x)$ is an iterated Ore extension. Proposition~\ref{prop:iteratedOre} thus shows that $x\in \Supp_{\mD_1}^{\hyp}(M)$ if and only if $x\in \Supp_{\mD_2}^{\hyp}(M)$.
\end{proof}

\subsection{Connection between homological and hypersurface support}

Fix an $n$-deformation datum $\mD$ with $Q\tto R$ as above. 

\subsubsection{}\label{sec:HH} By factoring out the ideal generated by $\fm^2$ in $Q$, we obtain a family of infinitesimal deformations of $R$, parametrised by $W^\ast$. This defines a map from $W^\ast$ to the second Hochschild cohomology group $\mathrm{HH}^2(R)$, see \cite{Ge}. Concretely, set 
\[J:=\fm Q/\fm^2 Q\subset Q/\fm^2 Q=:Q^{(1)}.\]
Choose a $\bk$-linear section $s:R\to Q^{(1)}$ of the algebra projection $Q^{(1)}\tto R$, with $s(1)=1$. Then we can define the $\bk$-linear morphism
\[\mathrm{def}_s: R\otimes R\to J,\quad a\otimes b\mapsto s(a)s(b)-s(ab).\]
By flatness of $Q$ over $Z$, we have a $\bk$-linear (actually $R$-bilinear) isomorphism
\[W\otimes R\xrightarrow{\sim} J,\quad v\otimes a\mapsto v \,s(a).\]
Note that the right-hand side does not depend on $s$. Composing the inverse of this isomorphism with $\mathrm{def}_s$ produces  a morphism $\mathrm{def}_s':R\otimes R\to W\otimes R$ which we can rewrite as a morphism
\[W^\ast\to\Hom_{\bk}(R\otimes R,R).\]
Using associativity of $Q$, this map takes values in cocycles, so it produces a morphism $W^\ast\to\mathrm{HH}^2(R)$ (which is canonical because different sections $s$ lead to the same maps above up to coboundary). It follows that we obtain a graded algebra morphism (with $W^\ast$ in degree 2)
\begin{equation}\label{eq:WRen}\Sym W^\ast\;\to\; \Ext^\bullet_{R^{\mathrm{en}}}(R,R),\end{equation}
with $R^{\mathrm{en}}=R\otimes R^{\mathrm{op}}$. This in turn leads to a central graded algebra morphism,
\begin{equation}\label{eq:algmor}\Sym W^\ast\;\to\; \Ext^\bullet_{R}(M,M),\end{equation}
for every $R$-module $M$, see \cite[\S 3.4]{NP2} or \cite[\S 4.3]{Lo}.


\subsubsection{}\label{sec:DHopf}
Now we assume that our deformation datum $\mD$ has some additional structure. We assume that $Q$ and $R$ are Hopf algebras, $Q\tto R$ is a morphism of Hopf algebras, and $Z\subset Q$ is a left coideal subalgebra. As a special case of \eqref{eq:algmor}, we have an algebra morphism
\[\Sym W^\ast\;\to\; E_R=\Ext^{2\bullet}_R(\bk,\bk).\]
By \cite[Theorem~4.8]{NP2} this is a finite morphism, which thus induces a morphism of varieties
\[\kappa:\mProj E_R\;\to\; \mP W.\]
To make everything explicit, we can trace the definition to observe that the algebra morphism is induced from the linear morphism
\begin{equation}\label{eq:explWH}
W^\ast\to H^2(R,\bk)=\Ext^2_R(\bk,\bk),\quad \lambda\mapsto [(\lambda\otimes \epsilon)\circ \mathrm{def}'_s].
\end{equation}
Here, we compute $H^\bullet(R,\bk)$ via the bar resolution, and $(\lambda\otimes \epsilon)\circ \mathrm{def}'_s$ in $\Hom_{\bk}(R\otimes R,\bk)$ only depends on $s$ up to coboundaries.

\begin{theorem}\cite[Theorem~7.1]{NP2}\label{thm:NP}
For any finite dimensional $R$-module $M$, we have
\[\Supp_{\mD}^{\hyp}(M)\;=\;\kappa\left(\cV_R(M)\right).\]
\end{theorem}

\begin{remark}\label{rem:nat}
Consider the set-up from \ref{sec:DD}. It follows from construction that for every $R_2$-module $M$, we obtain a commutative square with the morphisms from \eqref{eq:algmor}, the algebra morphism induced from $W_{n_2}^\ast\tto W_{n_1}^\ast$ and the restriction algebra morphism $\Ext^\bullet_{R_2}(M,M)\to\Ext^\bullet_{R_1}(M,M)$. Consequently, if we are further in the situation that both deformation data have the structure from \ref{sec:DHopf} and that $Q_1\to Q_2$ is a morphism of Hopf algebras, we obtain a commutative square
\[\xymatrix{\mProj E_{R_1}\ar[rr]^{\kappa_1}\ar[d]&& \mP W_{n_1}\ar[d]\\
\mProj E_{R_2}\ar[rr]^{\kappa_2}&& \mP W_{n_2} .}\]
For completeness we point out that the left vertical map is well-defined, since $E_{R_2}\to E_{R_1}$ is a finite algebra morphism, as follows from finiteness of the composed $\Sym W_{n_2}^\ast\to E_{R_1}$ and commutativity of the diagram.
\end{remark}

\begin{remark}\label{rem:grad}
Consider an $n$-deformation datum $\mD$ (given by $Z\to Q$) and assume that $Q$ has a $\mZ$-grading as a $\bk$-algebra, such that every $y_i\in Z\subset Q$ is homogeneous. Then $W=\fm/\fm^2$ is naturally a $\mZ$-graded vector space. Moreover, $R=Q/\fm Q$ is a $\mZ$-graded algebra.

Now assume also that $\mD$ has the additional structure from \ref{sec:DHopf} such that $Q$ and $R$ are graded Hopf algebras. Then also $\Ext^2_R(\bk,\bk)=H^2(R,\bk)$ inherits a $\mZ$-grading. It follows by choosing a graded section $s$ and analysing equation~\eqref{eq:explWH} that the morphism $W^\ast\to H^2(R,\bk)$ from \ref{sec:DHopf} is graded, where $W^\ast$ is understood to mean the (monoidal) dual of $W$ as a graded vector space.
\end{remark}


\section{The small quantum Borel}\label{sec:Borel}

In this section we work over $\bk=\mC$.

\subsection{Set-up}
\subsubsection{}\label{sec:data}
We fix $r\in\mZ_{>0}$ and a Cartan matrix $(a_{ij})\in \Mat_r(\mZ)$ of finite type. We denote the corresponding Lie algebra by $\fg$ with Borel subalgebra $\fb$ and nilpotent radical $\fn$. Set $n=\dim \fn$.

We let $\zeta\in\mC$ be a primitive $\ell$-th root of 1, where $\ell$ is odd, $\ell>h$ (with $h$ the maximum of the Coxeter numbers of the factors of $\fg$) and $\ell$ is prime to $3$ if $\fg$ has a factor of type $G_2$. Let $d_i$ be integers, with greatest common divisor 1 on each component, such that $(d_ia_{ij})$ is symmetric and positive definite and set $q_i:=\zeta^{d_i}$.

We also fix a non-empty subset $I\subset \{1,\ldots, r\}$, leading to a (semisimple Levi) subalgebra $\fg_I$ generated by the (positive and negative) root vectors corresponding to $I$, with Borel subalgebra $\fb_I=\fb\cap \fg_I$ with nilpotent radical $\fn_I=\fn\cap\fg_I$. Set $m=\dim \fn_I$. Note that $\ell$ inherits our assumptions (from the previous paragraph) for $\fg_I$ (although keeping $q_i$ as before might lead to the consideration of a different primitive root $\zeta'$ of 1 when the $d_i$ have non-trivial common divisor on some component).

\subsubsection{}
We denote by $U^{\DK}_\zeta(\fg)$ the De Concini-Kac form of the quantum group, generated as a $\mC$-algebra by $E_i, F_i, K_i^{\pm 1}$, for $1\le i\le r$, subject to relations \cite[1.2.1-1.2.5]{DK}. Its (Hopf) subalgebra $U^{\DK}_\zeta(\fb)$ is generated by $E_i,K_i^{\pm 1}$.

We choose a reduced expression of the longest element of the Weyl group of $\fg$ that begins with a reduced expression of the longest element of the Weyl group of $\fg_I$. As in \cite[\S1]{DK} this gives `root vectors' $E(j)\in U_\zeta^{\DK}(\fb)$, for $1\le j\le n$, where $E(j),1\le j\le m$ are the corresponding root vectors for $U_\zeta^{\DK}(\fb_I)$. Similarly we have negative root vectors $F(j)$.
By \cite[Lemma~3.1]{DK}, the elements $E(j)^{\ell}$, $F(j)^\ell$ and $K_i^{\ell}$ are central in $U^{\DK}_\zeta(\fg)$.

We denote the quotient of $U^{\DK}_\zeta(\fg)$ with respect to $E(j)^{\ell}$, $K_i^{\ell}-1$ and $F(j)^\ell$ by $u_\zeta(\fg)$. This is Lusztig's small quantum group, which is also the image of the canonical Hopf algebra morphism $U_\zeta^{\DK}(\fg)\to U_\zeta^{\L}(\fg)$, where $U_\zeta^{\L}(\fg)$ is Lusztig's (divided power) form of the big quantum group from \cite{Lu} (with $K_i^\ell=1$ imposed). Similarly, we have the quotient $u_\zeta(\fb)$ of $U^{\DK}_\zeta(\fb)$ by the relations  $E(j)^{\ell}$ and  $K_i^{\ell}-1$.

With respect to $I$, we have (Hopf) inclusions of all quantum groups\footnote{There are some subtleties here. Firstly, $\ell$ does inherit our assumptions from \ref{sec:data} with respect to the Cartan matrix $(a_{ij})_{i,j\in I}$. However, it can happen that the $d_i$ have a non-trivial common divisor on one of the components (though by construction such a divisor will always be coprime to $\ell$). Consequently, the notation $u_\zeta(\fg_I)$ is slightly misleading as the Hopf algebra is actually of the form $\otimes_s u_{\zeta^{(s)}}(\fg^{(s)})$, where the $\fg^{(s)}$ are the simple factors of $\fg_I$ and every $\zeta^{(s)}$ is some primitive $\ell$-th root of $1$, not necessarily equal to $\zeta$.}, for instance $u_\zeta(\fb_I)\subset u_\zeta(\fb)$ and $U^{\DK}_\zeta(\fg_I)\subset U^{\DK}_\zeta(\fg)$.

\subsubsection{}\label{sec:cohomb} It is proved in \cite[Theorem~2.5]{GK} that
\[E_{u_\zeta(\fb)}=\Ext^{2\bullet}_{u_\zeta(\fb)}(\mC,\mC)=\Ext^{\bullet}_{u_\zeta(\fb)}(\mC,\mC)\cong\mC[\fn]=\Sym(\fn^\ast),\]
where the image of $\fn^\ast$ in $E_{u_\zeta(\fb)}$ lives in homological degree 2. In particular $\mProj E_{u_\zeta(\fb)}\cong \mP\fn$. Moreover, the restriction morphism $E_{u_\zeta(\fb)}\to E_{u_\zeta(\fb_I)}$ corresponds to $\fn^\ast\tto \fn^\ast_I$.

\subsection{Connection to Ore extensions and deformation data}

\subsubsection{}\label{def:Qb} We set $Q:=U_{\zeta}^{\DK}(\fb)/(K_i^{\ell}-1)$, and define the central algebra morphism
\[Z_n=\bk[y_1,\ldots, y_n]\;\to\; Q,\quad y_j\mapsto E(j)^{\ell}.\]
We also let $Q_I$ be the Hopf subalgebra of $Q$ generated by $E_i,i\in I$ (equivalently $E(j)$ for $j\le m$), and \emph{all} $K_i, 1\le i\le r$.
Similarly, we define the central algebra morphism
\[Z_m=\bk[y_1,\ldots, y_m]\;\to\; Q_I,\quad y_j\mapsto E(j)^{\ell}.\]

The following lemma shows that these algebra morphisms are deformation data, which clearly fit into the set-up of \ref{sec:DD}.
\begin{lemma}
\begin{enumerate}
\item $\mC[K_1,\ldots, K_{r}]/(K_i^{\ell}-1)\hookrightarrow Q_I$ and $Q_I\hookrightarrow Q$ are iterated Ore extensions.
\item $Q$ and $Q_I$ have finite global dimension.
\item $Q$ is finite free over $Z_n$ and $Q_I$ is finite free over $Z_m$.
\end{enumerate}
\end{lemma}
\begin{proof}
Part (1) follows from \cite[\S 1.7]{DK}, since a chain of Ore extensions $Q(0)\subset Q(1)\subset \cdots\subset Q(n)=Q$ passing through $Q_I=Q(m)$ is given by the subalgebras
\[Q(j)=\langle K_1,\cdots, K_r, E(1),\ldots E(j)\rangle \subset Q. \]

Part (2) follows from Proposition~\ref{prop:iteratedOre} and part (1), since the source algebra in (1) is semisimple.

Part (3) follows from \cite[Corollary~3.3]{DK}.
\end{proof}
\subsubsection{} We have $Q/(\fm_n)\cong u_\zeta(\fb)$ canonically and $Q$ inherits a Hopf algebra structure from $U_\zeta^{\DK}(\fb)$ such that $Q\tto u_\zeta(\fb)$ is a Hopf algebra morphism. Moreover, $Z_n$ and $Z_m$ are left coideals, even Hopf subalgebras (see \cite[Proposition~5.6]{DKP}), so the deformation system satisfies the extra structure from \ref{sec:DHopf}.

Similarly, the Hopf algebra
\[Q_I/(\fm_m)\;\cong\; u_\zeta(\fb_I)\mathrel{\ooalign{$\otimes$\cr$\rtimes$}} \mC[K_i\mid i\not\in I]/(K_i^\ell-1)\]
has $u_\zeta(\fb_I)$ as a normal Hopf subalgebra. We can think of the right-hand side as $u_\zeta(\tilde{\fb}_I)$, where $\tilde{\fb}_I\subset\fg$ is the Lie subalgebra spanned by $\fb_I$ and the Cartan subalgebra of $\fg$. In order to apply Lemma~\ref{lem:Hopfnormal}, we require the following lemma.
\begin{lemma}\label{lem:trivact}
The action of $\mC[K_i\mid i\not\in I]/(K_i^\ell-1)$ on $E_{u_\zeta(\fb_I)}$ is trivial.
\end{lemma}
\begin{proof}
The approach of \cite{GK} is to compute $H^\bullet(u_\zeta(\fb_I),\mC)$ as the  $\mC[K_i\mid i\in I]/(K_i^\ell-1)$-invariants in $H^\bullet(u_\zeta(\fn_I),\mC)$. It is clear from the construction {\it loc. cit.} that $\mC[K_i\mid 1\le i\le r]/(K_i^\ell-1)$ acts trivially on the $\mC[K_i\mid i\in I]/(K_i^\ell-1)$-invariants. Equivalently, the surjection $E_{u_\zeta(\fb)}\tto E_{u_\zeta(\fb_I)}$ from \ref{sec:cohomb} is equivariant for the action of the $K_i$'s, but their action on the left-hand side is manifestly trivial.
\end{proof}

\subsection{Naturality of cohomological support}

\begin{theorem}\label{thm:B}
For every finite dimensional $u_\zeta(\fb)$-module $M$,
\[\cV_{u_\zeta(\fb_I)}(M)\;=\;\cV_{u_\zeta(\fb)}(M)\cap \mP\fn_I.\]
\end{theorem}
\begin{proof}
By Lemmas \ref{lem:Hopfnormal} and \ref{lem:trivact}, we can instead prove
\[\cV_{u_\zeta(\tilde\fb_I)}(M)\;=\;\cV_{u_\zeta(\fb)}(M)\cap \mP\fn_I.\]
The latter follows from Proposition~\ref{prop:iteratedOre2} and the naturality in Remark~\ref{rem:nat}, since in our set-up $\kappa_1,\kappa_2$ are bijections. Indeed, the corresponding graded algebra morphisms
\[\Sym W_n^\ast\to E_{u_\zeta(\fb)}=\Sym \fn^\ast \quad\mbox{and}\quad \Sym W_m^\ast\to E_{u_\zeta(\tilde\fb_I)}=\Sym \fn_I^\ast  \]
from \ref{sec:DHopf} must be isomorphisms, since $W_n\cong \fn$ and $W_m\cong \fn_I$ as vector spaces, and the algebra morphisms are finite.
\end{proof}


\section{Naturality of support for tilting modules}\label{sec:Full}
In this section, we keep the data, notation and assumptions from \ref{sec:data}.

\subsection{Reducing to the Borel subalgebra}

\subsubsection{} It was proved in \cite{GK} that
\[E_{u_\zeta(\fg)}=\Ext^{2\bullet}_{u_\zeta(\fg)}(\mC,\mC)=\Ext^{\bullet}_{u_\zeta(\fg)}(\mC,\mC)\cong\mC[\cN],\]
where $\cN\subset\fg$ is the nilpotent cone. The algebra isomorphism is graded, with degree on $\mC[\cN]$ coming from the conical nature of $\cN$, up to a factor $2$. Let $G$ be the connected algebraic group of adjoint type, acting faithfully on~$\fg$. 

To link more directly with the literature, we use the affine version of the cohomological support. Concretely, for a non-zero $u_\zeta(\fg)$-module $M$, we will write $\widehat{\cV}_{u_\zeta(\fg)}(M)\subset \cN$ for the cone corresponding to ${\cV}_{u_\zeta(\fg)}(M)$, while for the zero module $M$, we set $\widehat{\cV}_{u_\zeta(\fg)}(M)=\varnothing$. This corresponds to the usual definition based more directly on $\operatorname{mSpec}E_{u_\zeta(\fg)}$. We do the same for $u_\zeta(\fb)$.

\begin{theorem}\label{thm:gb}
Let $M$ be a finite-dimensional type 1 module of $U_\zeta^{\L}(\fg)$. Then
\[\widehat\cV_{u_\zeta(\fb)}(M)\;=\;\widehat\cV_{u_\zeta(\fg)}(M)\cap \fn.\]
 \end{theorem}

\begin{remark}
Theorem~\ref{thm:gb} is known, see \cite[Theorem~7.3.3]{BKN}. However, since the proof in the current arxiv version (v3) contains a minor gap to be fixed later, we provide an independent proof. The big outline follows the approach from \cite{BKN}, but inserts some ideas from hypersurface support and Ore extensions to fill in the details.
\end{remark}

 \subsubsection{}Choose $1\le i\le r$ and consider the corresponding minimal parabolic subalgebra $\fb\subset\ffp\subset\fg$. We consider $Q':= U^{\DK}_\zeta(\ffp)/(K_j^\ell-1,1\le j\le r)$. It follows quickly from \cite[(1.2.2), (1.2.3)]{DK} that $Q\subset Q'$ is an Ore extension, with $Q$ from \ref{def:Qb}. Moreover, the extension fits into a commutative diagram
 \[\xymatrix{
Z_n\ar[rr]\ar[d]&& Q\ar[d]\\
Z_{n+1} \ar[rr]&& Q',
 }\]
 where the bottom arrow sends $y_{n+1}$ to $F_i^{\ell}$, and corresponds to a deformation datum $\mD'$. The small quantum group $u_\zeta(\ffp)$ corresponding to $\ffp$ is the finite dimensional algebra in this deformation datum. Moreover, since $F_i^{\ell}$ is primitive in $Q'$, the deformation datum has the additional structure from \ref{sec:DHopf}. As follows from \cite{GK, BNPP}, we have $H^{2\bullet}(u_\zeta(\ffp),\mC)\cong \mC[\cN\cap\ffp]$ (equivariantly under the actions of the adjoint group for $\ffp$), where the restriction morphism to $H^{2\bullet}(u_\zeta(\fb),\mC)$ is simply restriction of functions to $\fn$.

 \begin{proposition}\label{prop:pi}
 For a finite dimensional $u_\zeta(\ffp)$-module $M$, we have
 \[\cV_{u_\zeta(\fb)}(M)\;=\; \cV_{u_\zeta(\ffp)}(M)\cap\mP\fn.\]
 \end{proposition}
\begin{proof}
By a direct application of Proposition~\ref{prop:iteratedOre2}, we have
\[\Supp^{\hyp}_{\mD}(M)\;=\;\Supp^{\hyp}_{\mD'}(M)\cap \mP W_{n}.\]
Moreover, the commutative diagram from Remark~\ref{rem:nat} is
\[\xymatrix{\mP\fn\ar[rr]^{\sim}\ar@{^{(}->}[d]&&\mP W_n\ar@{^{(}->}[d]\\
P\ar[rr]^{\kappa}&&\mP W_{n+1},}\]
where $P$ is the projectivisation of the cone $\cN\cap\ffp$. Theorem~\ref{thm:NP} shows that $\kappa(\cV_{u_\zeta(\ffp)}(M))\cap \mP W_n$ equals the image of $\cV_{u_\zeta(\fb)}(M)$ under the top isomorphism. The proposition thus follows if we prove the following claim: only elements of $\mP\fn\subset P$ are sent to $\mP W_n\subset\mP W_{n+1}$ under $\kappa$.

Now we prove this claim. Rather than working with the above diagram of varieties, we can work with the diagram of linear morphisms corresponding to the degree 2 component of the defining algebra morphisms
\[\xymatrix{\fn\ar[rr]^{\sim}\ar@{^{(}->}[d]&&W_n\ar@{^{(}->}[d]\\
\ffp'\ar[rr]^{\kappa}&&W_{n+1},}\]
where the linear map is again denoted by $\kappa$ with slight abuse of notation, and $\ffp'=[\ffp,\ffp]$. By dimension count, $\kappa$ is not injective, but by finiteness of the algebra morphism it is surjective, so $\kappa$ has a one-dimensional kernel $L\subset\ffp'$.
Consider the standard $\mathfrak{sl}_2$-triple $(e, h, f)$ in $\ffp$ where $e$ resp. $f$ is the positive resp. negative root vector for the root corresponding to $i$. Then we can realise
\[\cN\cap\ffp\;=\;\{v+ae+bh+cf\mid b^2+ac=0\}\;\subset\; \ffp'\;=\; \{v+ae+bh+cf\},\]
with $v$ in the nilpotent radical of $\ffp$ and $a,b,c\in\mC$. In particular, what we need to prove in order to arrive at the claim is that $L\subset \fn\oplus \mC h$.

To prove $L\subset \fn\oplus \mC h$, we will set up a graded version of $\kappa$. We consider $Q'$ as a $\mZ$-graded algebra where the degree of $E_j$ is $a_{ij}\in\mZ$, the degree of $F_i$ is $-2=-a_{ii}$ and the degree of $K_j$ is $0$, for all $1\le j\le r$. This makes $Q$ a graded subalgebra, and the generators of $Z_{n+1}$ are homogeneous. It follows via Remark~\ref{rem:grad} that the above diagram can be viewed in the category of graded vector spaces. We thus have resulting gradings on $W_n$ and $W_{n+1}$ where the non-zero homogeneous parts of $W_n$ live in degrees that are sums of integers $\ell a_{ij}$, with $1\le j\le r$, and $W_{n+1}$ is the direct sum of $W_n$ with a one-dimensional space in degree $-2\ell$ (corresponding to $y_{n+1}=F_i^{\ell}$). To prove $L\subset \fn\oplus \mC h$ it is thus sufficient to show that $h\in\ffp'$ lives in degree $0$, as that then implies $\kappa(h)\in W_n=\kappa(\fn)$. 

Finally, we prove that $h\in\ffp'$ lives in degree $0$. By construction, the grading on $\ffp'$ originates from
\begin{equation}\label{eq:pH2}\ffp'\;\cong\; \left(H^2(u_\zeta(\ffp),\mC)\right)^\ast\end{equation}
and the grading on $u_\zeta(\ffp)$ described in the previous paragraph. It is thus actually sufficient to show that this grading on $H^2(u_\zeta(\ffp),\mC)$ equals, up to a multiplicative shift by $\ell$, the grading corresponding to the action of the $SL_2$-torus $T\cong \mG_m$ (with $\operatorname{Lie}T=\mC h$).
Indeed, this action on the right-hand side of \eqref{eq:pH2} (for the entire adjoint group for $\ffp$) is transferred by the isomorphism to the usual adjoint action on $\ffp'$ (so that $T$ acts trivially on $h\in\ffp'$). Recall from \cite[Chapter~2]{Dr1} that the adjoint group action on cohomology of the small quantum group is obtained from the (right) adjoint action of $U^{\L}_\zeta(\ffp)$ on $u_\zeta(\ffp)$, which on cohomology descends to an action of $U^{\L}_\zeta(\ffp)/\hspace{-1mm}/u_\zeta(\ffp)\cong U(\ffp)$, the ordinary enveloping algebra of $\ffp$. We must thus find an element $H\in U^{\L}_\zeta(\ffp)$ that is sent to $h\in\operatorname{Lie} T$ under $U^{\L}_\zeta(\ffp)\tto U(\ffp)$ and for which the right adjoint action of $H$ on $u_\zeta(\ffp)$ sends $x\in u_\zeta(\ffp)_{d}$ (a homogeneous element in degree $d$) to $- d x/\ell$. Such an element is constructed in \cite[Theorem~3.1(d)]{Le}. More explicitly, it is a primitive element
\[H\;=\;\left[\begin{smallmatrix} K_i;0 \\ \ell \end{smallmatrix}\right]+f(K_i^2)\;\in\; U^{\L}_\zeta(\ffp),\]
where $f\in\mathbb{Q}[t]$ is a polynomial with $f(1)=0$, and its right adjoint action defines the grading on $Q'$ as desired.
\end{proof}

 \begin{remark}
 Unfortunately, we cannot continue the procedure of Proposition~\ref{prop:pi} to arrive at a version of Theorem~\ref{thm:gb} that does not rely on a lift to $U^{\L}_\zeta(\fg)$, since the only obvious path to continue with Ore extensions $Q\subset Q'\subset\cdots$ corresponds to adding $F_j$ with $a_{ij}=0$.
 \end{remark}
 
 \begin{proof}[Proof of Theorem~\ref{thm:gb}]
 We can start from the weaker version of the claim, as proved in \cite[Theorem~6.1]{Dr}:
 \[\widehat\cV_{u_\zeta(\fg)}(M)\;=\; G\cdot\widehat\cV_{u_\zeta(\fb)}(M).\]
 For the suggested equality in the theorem, the left-hand side is included in the right-hand side via the weak naturality that follows immediately from the definitions. Thus we consider \[x\in \widehat\cV_{u_\zeta(\fg)}(M)\cap \fn\;=\; G\cdot\widehat\cV_{u_\zeta(\fb)}(M) \,\cap\, \fn.\] Choose $g\in G$ and $y\in \widehat\cV_{u_\zeta(\fb)}(M)$, to write $x=g\cdot y$. Use the Bruhat decomposition to write $g=b_1\dot{w}b_2$, for $b_1,b_2\in B$ and $w$ in the Weyl group. By $B$-invariance of $\widehat\cV_{u_\zeta(\fb)}(M)$, see \cite{Dr1, GK}, we can rewrite things with $b_2=1$ and conclude that $\dot{w}\cdot y\in\fn$.

 If $w$ is a simple reflection, we write $w=s_i$ for some $1\le i\le r$. Since $\cV_{u_\zeta(\ffp)}(M)$ is invariant under the action of $\dot{s}_i$, the conclusion of the previous paragraph in combination with Proposition~\ref{prop:pi} shows that $\dot{w}\cdot y\in \widehat\cV_{u_\zeta(\fb)}(M)$. The case of general $w$ can then be proved by a clean inductive argument, as explained on p29 of \cite{BKN}.
 \end{proof}

\subsection{Main results}

It is well-known that for a $U_\zeta^{\L}$-module, the subset $\widehat\cV_{u_\zeta(\fg)}(M)\subset \cN$ is closed under the $G$-action, see \cite{GK, Dr1}. It is also known that $E_{u_\zeta(\fg)}\to E_{u_\zeta(\fg_I)}$ corresponds to the canonical inclusion $\cN(\fg_I)\hookrightarrow\cN$, see \cite[\S 4.3.1]{CEO}. In fact, by further restricting to $u_\zeta(\fb_I)$ and using $G_I$-equivariance it can be reduced to the Borel case. Here, we let $G_I<G$ be the subgroup generated by the root subgroups corresponding to $\fg_I$. It might not be of adjoint type.

\begin{theorem}\label{thm:G}
Let $M$ be a finite-dimensional type 1 module of $U_\zeta^{\L}(\fg)$, viewed as a $u_\zeta(\fg)$-module. Then
\[\widehat\cV_{u_\zeta(\fg_I)}(M)\;=\;\widehat\cV_{u_\zeta(\fg)}(M)\cap\fg_I.\]
\end{theorem}
\begin{proof}
By Theorem~\ref{thm:gb}
\[\widehat\cV_{u_\zeta(\fb)}(M)\;=\;\widehat\cV_{u_\zeta(\fg)}(M)\cap \fn\quad\mbox{and}\quad \widehat\cV_{u_\zeta(\fb_I)}(M)\;=\; \widehat\cV_{u_\zeta(\fg_I)}(M)\cap\fn_I.\] 
The equality now follows from Theorem~\ref{thm:B}, since two $G_I$-invariant subsets of $\cN(\fg_I)$ are equal if and only if their intersections with $\fn_I$ are equal (by the principle $\cN=G\cdot\fn$).
\end{proof}

Let $\Tilt U_\zeta(\fg)$ be the monoidal category of type 1 tilting modules of $U_\zeta(\fg)$, see \cite{An}. Tensor ideals of abelian type, see \cite[Definition~2.1.1]{CEO}, in $\Tilt U_\zeta(\fg)$ were classified in  \cite[Theorem~5.4.3]{CEO} conditionally on the validity of the naturality in support now proved.

\begin{theorem}
\begin{enumerate}
\item The equivalent conditions in \cite[Proposition~5.3.1]{CEO} are all true in the generality stated there, equivalently \cite[Conjecture~5.3.3]{CEO} is true.
\item The tensor ideals of abelian type in $\Tilt U_\zeta^{\L}(\fg)$ are in bijection with nilpotent $G$-orbits in the nilpotent cone of $\fg$.
\end{enumerate}
\end{theorem}
\begin{proof}
The second of the equivalent conditions in \cite[Proposition~5.3.1]{CEO} is precisely what is proved in Theorem~\ref{thm:G}.
Part (2) was proved in \cite[Theorem~5.4.3]{CEO}, conditional on the validity of the now verified \cite[Conjecture~5.3.3]{CEO}.
\end{proof}


\appendix

\section{Connection to the formal setting}\label{sec:app}

In the main text we did not work in the formal set-up from \cite{NP1, NP2}. To justify our usage of those results we thus require two facts. Firstly, our notion of deformation datum from \ref{def:datum} can be completed to the set-up from \cite[\S 2.1- 2.2]{NP1}. Secondly, the resulting hypersurface support for the completed version is identical to the one we use, including the comparison map in \eqref{eq:WRen}.

For the remainder of the appendix, fix an algebraically closed field $\bk$ and an $n$-deformation datum~$\mD$ as in \ref{def:datum}, corresponding to~$Q,R,\fm$.

\subsubsection{}
Consider the complete local algebra $\widehat{Z}=\bk[[y_1,\ldots, y_n]]$ with maximal ideal $\widehat{\fm}$, and set
\[\widehat{Q}=\widehat{Z}\otimes_{Z}Q\cong \varprojlim Q/\fm^i Q.\]
We observe that $\widehat{Z}\to \widehat{Q}$ satisfies all assumptions required from \cite[\S 2.1- 2.2]{NP1}. Firstly, finite freeness is inherited and we have $\widehat{Q}/\widehat{\fm}\widehat{Q}\cong R$ canonically. Again, $\widehat{Q}$ is noetherian, as a finite algebra over the commutative noetherian ring $\widehat{Z}$. The semisimple quotient of $R$ is separable since $\bk$ is algebraically closed. The main condition to verify is that $\widehat{Q}$ still has finite global dimension (so in particular is Gorenstein). The special case directly relevant to the current paper is actually proved in \cite[Lemma~2.10]{NP2} and the proof applies to our generality from~\ref{def:datum}. Finally, the required `integral form' of $\widehat{Z}\to\widehat Q$ is precisely the original $Z\to Q$.

\subsubsection{}
We have canonically
\[\widehat{W}:=\widehat{\fm}/\widehat{\fm}^2\cong \fm/\fm^2 = W.\]
Hypersurface support in \cite{NP1, NP2} is defined precisely as we outlined in \S\ref{sec:hyper}, but with respect to the lift to $\widehat{Q}$ and taking values in $\mP\widehat{W}$.
The following lemma shows this does not make a difference, so that in particular, the definition in \ref{def:hyper} makes sense.
\begin{lemma}
Let $M$ be a finite-dimensional $R$-module and let $f\in\fm\subset\widehat\fm$ have non-zero class in $\fm/\fm^2$, then $\pd_{\widehat{Q}/(f)}M=\pd_{Q/(f)}M.$
\end{lemma}
\begin{proof} 
Set $C=Z/(f)$ with $\fq\subset C$ the image of $\fm\subset Z$, and $\widehat{C}=\widehat{Z}/(f)$. All localisation will be with respect to $C$, so for any $\fp\in\Spec C$, we have $(-)_{\fp}=C_{\fp}\otimes_C-$. For example $M_{\fq}=M$ and  $M_{\fp}=0$ whenever $\fq\not=\fp\in\Spec C$, since $\fq$ acts trivially on $M$ by construction.
We prove the equality in the lemma in two steps, by considering the intermediate algebra $Q/(f)\to A\subset \widehat{Q}/(f)$, with $A:= (Q/(f))_{\fq}$.


First we show that $\pd_{Q/(f)}M=\pd_A M$. Since $Q/(f)$ is noetherian, we can consider a projective resolution of $M$ that is finitely generated in every degree, which leads to 
\[\Ext^\bullet_{Q/(f)}(M,-)_{\fp}\cong \Ext^\bullet_{Q/(f)}(M,(-)_{\fp})\cong\begin{cases}
0 &\mbox{ if } \fp\not=\fq,\\
\Ext^\bullet_{A}(M,(-)_{\fq}) &\mbox{ if } \fp=\fq.
\end{cases}\]
This implies the desired equality.

The equality $\pd_{\widehat Q/(f)}M=\pd_A M$ follows similarly, since $\widehat Q/(f)\cong A\otimes_{C_{\fq}}\widehat{C}$, where $C_{\fq}\to \widehat{C}$ is faithfully flat.
\end{proof}

\subsubsection{} Finally, we need to verify that the defining degree 2 morphism for \eqref{eq:WRen} is the same if we had been working with $\widehat{Z}\to\widehat{Q}$ instead. This is straightforward, since $\widehat{Q}/ \widehat{\fm}^2\widehat{Q}\cong  Q/\fm^2 Q$ canonically.

\subsection*{Acknowledgements}
 The research was partly supported by the ARC grants FT220100125 and DP250100762. The author thanks Daniel Nakano for pointing out that Theorem~\ref{thm:gb} reduces \cite[Conjecture~5.3.3]{CEO} to Theorem~\ref{thm:B}.

\end{document}